\documentclass{engine/svjour3}
\smartqed  
\usepackage{epsfig}
\usepackage{amsmath}
\usepackage{amssymb}
\usepackage{mathrsfs}
\usepackage{chemarrow}
\usepackage[all]{xy}
\usepackage{times}
\usepackage{comment}
\usepackage{graphicx}
\usepackage{lmodern}
\usepackage[utf8]{inputenc}
\usepackage[T1]{fontenc}
\usepackage{textcomp}
\usepackage{bm}
\usepackage[hyphens]{url}
\usepackage{breakurl}

\usepackage{hyperref}
\hypersetup{
    colorlinks=true,
    linkcolor=blue,
    filecolor=magenta,      
    urlcolor=cyan,
    pdftitle={Sharelatex Example},
}

\def\v1{\mbox{\boldmath$\big|1\big\rangle$}}

\usepackage[
backend=biber,
style=alphabetic,
sorting=ynt
]{biblatex}
\begin{document}

\title{A List of Open Problems in Number Theory Posed by Ibn al-Khawwām (13th Century): Historical and Arithmetic Analysis}
\titlerunning{Unsolved arithmetic problems of the time of Ibn Khawwām} 

\author{K.I.A Derouiche
}

\institute{K.I.A Derouiche \at
BMW Motorrad Algeria, IT.Departement\\
 8, Rue Mohammed Loubi,\\
Alger, Algérie \\ 
\email{kamel.derouiche@gmail.com \\
	k.derouiche@algeriemotors.com
}}

\dedication{DRAFT PAPER}

\maketitle

%%%%%%%%%%%%%%%%%

\begin{abstract}
Mathematicians have long been fascinated by the resolution of algebraic and Diophantine 
equations in search of integer or rational solutions. This article presents a list of thirty-three open 
problems in number theory, posed in the 13th century by Ibn al-Khawwām al-Baghdādī (Abdallāh 
ibn Muḥammad ibn Muḥammad al-Khawwām), extracted from his arithmetic treatise al-Fawāid 
al-Bahāiyya fī qawāid al-ḥisābiyya. We provide a historical and arithmetic analysis of these 
problems and their solutions, and, whenever possible, offer original solutions and bibliographic 
notes. This work situates these problems within the broader development of number theory and 
highlights their continuing mathematical relevance
\end{abstract}


\tableofcontents

%%%%%%%%%%%%%%%%%%%%%%%

\section{Introduction}
The term open problem or unsolved problem in mathematics generally refers to a question or 
conjecture that has remained unresolved for a long period of time \cite{KOSYVAS2010}.
A classical example is the Poincaré Conjecture, stated in 1904 in his Fifth Complement to Analysis 
Situs — “Is every compact, boundaryless, simply connected 3-manifold homeomorphic to the 
3-sphere?” \cite{Poincare1904} — which was solved a century later by Grigori Perelman 
\cite{Perelman2003}. Another emblematic example is Fermat’s Equation, $x^{n} + y^{n} = z^{n}$ for 
$n > 3$, proved by Andrew Wiles after more than three centuries of research \cite{Wiles1995, 
Taylor1995}.

Likewise, the Riemann Hypothesis, concerning the distribution of prime numbers, remains one of the major unsolved challenges in mathematics \cite{Riemann1859, Bombieri2022}.
These conjectures belong to the category of great open problems, which encompass all areas of mathematics and have often been compiled into famous lists — such as the one presented by David Hilbert in 1900 at the International Congress of Mathematicians in Paris \cite{Hilbert1902}.

In the twentieth century, Steve Smale extended this tradition by proposing eighteen open problems 
covering an even broader range of topics, including theoretical computer science, mathematical 
economics, partial differential equations, fluid mechanics, algebraic geometry, and even artificial 
intelligence \cite{Smale1998}.
Later, Misha Gromov formulated another set of open problems focusing on hyperbolic geometry, 
infinite groups, and metric geometry \cite{Gromov2014}.
Other, more specialized lists have also appeared, particularly in the fields of nonlinear analysis and 
optimization \cite{Hiriart-Urruty2007}.

In the same spirit, Jean-Louis Colliot-Thélène compiled a modern collection of open problems in 
algebraic geometry and number theory, often accompanied by detailed commentary and 
discussion \cite{Colliot2022}. These initiatives illustrate the enduring mathematical and 
philosophical interest in unsolved problems — genuine milestones in the evolution of scientific 
thought.

Within this long-standing tradition, an earlier contribution is devoted to the list of thirty-three 
arithmetic problems formulated by the thirteenth-century mathematician Ibn al-Khawwām 
al-Baghdādī. This list appears in his treatise al-Fawāid al-Bahāiyya fī qawāid al-ḥisābiyya ("The Bahāī 
Benefits on the Principles of Arithmetic"). A first study of these problems was conducted by Mahdi 
Abdeljouad and al\cite{Abdeljouad1986}, which highlighted the historical and 
scientific significance of these questions during the height of Islamic civilization and their 
connection to number theory.

The present paper continues this line of research. It offers both a historical analysis and an 
arithmetical examination of Ibn al-Khawwām’s thirty-three problems. We reproduce the original list 
and, whenever possible, provide analytical solutions and additional bibliographical notes.
Our aim is to place these problems within their historical context while emphasizing their 
importance in the development of number theory and their relevance to contemporary 
mathematical research.

\subsection{Organization of the Article}
The article is organized as follows: the second section is devoted to the biography and bibliographical elements concerning Ibn al-Khaww$\bar{a}$m, followed by a section presenting the various commentators on the list of thirty-three problems Ibn al-Khaww$\bar{a}$m. The fourth section discusses and sets out the thirty-three arithmetic problems.

\section{Bibliographical elements concerning Ibn al-Khaww\={a}m}
Abdallah ibn Muhhammad ibn Abd al-Razz\={a}q Al-kharab\={a}wi Ima\={a} El-Edine
Ibn al-Khaww\={a}m  born in 643H/1245 in Baghdad, was a physician, writer, philosopher, and
mathematician. He was nicknamed "the Arithmetician" in Baghdad by his peers, due to his
fame and his influence on the arithmetic of his era, which he maintained throughout his life, as
reported by the historian Shams ad-Din adh-Dhahabi \cite{Dhahabi2016}.

He was well-versed in the religious sciences, as was customary at the time.\footnote{At that
time, and even long before, education was conducted in higher colleges (madrasas) where
both religious subjects and secular disciplines (literature, grammar, astronomy, mathematics)
were taught. These colleges hosted professors and students for teaching activities on a wide
variety of topics \cite{Djebarr2012}.}

Because of his fame, Ibn al-Khaww\={a}m was widely sought after at court as a tutor for the
sons of ministers and princes in Isfahan, under the reign of Abaka Khan (d. 680H/1282), son of
H\"{u}lag\"{u} Khan. Around 475H/1277, at the age of thirty-two, in this same city, Ibn 
al-Khaww\={a}m wrote his book al-Faw\={a}'id al-Bah\={a}'iya fi qaw\={a}'id al-His\={a}biya 
during the month of Sha'bān (January 675H/1276) \cite{Fazlioglu1995}.

He taught Shafi'ite jurisprudence (fiqh) at Dar ad-Dhahab in Baghdad, directed the Mashyakh\={a}t 
ar-Riba of Baghdad, and remained active there until his death in 724H/1324. Among his students 
was the eminent mathematician and physicist Kamāl al-D\={a}n al-Fārisī (1267-1319) 
\cite{El-Fuwti1996}.

\subsection{The works of Ibn al-Khaww\={a}m}
The historians Shams ad-D\={i}n adh-Dhahabi and Ibn Hajar al-'Asqal\={a}ni mention only four works (\cite{Dhahabi2016}, \cite{Asqalani1972}) attributed to Ibn al-Khaww\={a}m:
\begin{enumerate}
	\item An Introduction to Medicine.
	\item AL-Tadkira Al-Sa\={a}diya Fi EL-KAWANIN AL-TYBIYA. This is a manuscript dealing with medicine, catalogued in the al-Awqāf library in Mosul. However, no other author has mentioned this work. We believe that this book could possibly be a version of his Introduction to Medicine; however, in our opinion, it does not appear to be the same work.
	\item Rissa\={a}lat Fi Al-Firassa\={a}(Epistle on Physiognomy).
	\item Ris\={a}la fi f\={a}m al-Maq\={a}la al-Ashir\={a} min Kit\={a}b Uql$\bar{u}$dis.
	\item Al-Faw\={a}'id al-Bah\={a}'iya fi qaw\={a}'id al-His\={a}biya Our present book).
\end{enumerate}

\section{The Commentators of d'AL-FAWA'ID AL-BAHA'IYA }
Ibn Khaww\={a}m \`a  book has attracted considerable interest from its appearance until today. Several commentators have addressed the fifth and final volume, which is devoted to the list of thirty-three problems some have treated only certain parts of it, others have commented mainly on questions of an algebraic nature, while yet others have merely suggested possible approaches. Another category of commentators reproduced the solutions that were provided and offered explanations, often with a pedagogical purpose. To the best of our knowledge, no commentator has actually solved the problems posed in this list, with the exception of the nearly complete work of Abdeljouad and Hemida \cite{Abdeljouad1986}, prior to the modern resolution.
\begin{itemize}
	\item Ibn al-Khaww\={a} Kam\={a}l ad-D\={i}n Hassan al-F\={a}risi(12 January  1267- 730H/1319), in his work entitled "Ass\={a}s al-qaw\={a}'id fi usul al-faw\={a}'id" (The Foundation of Rules in the Principles of Applications) \cite{Djebarr2006}, presents an introduction and five chapters dealing with arithmetic rules, notarial and commercial transactions, the measurement of areas and volumes, and the last two treatises focus on algebra.
	
	\item Abd Ali ibn Muhammad ibn Husayn Birjandi, the date of his death is not precisely known and deserves further clarification\cite{Kusuba2014}. Hajji Khalifa mentions three different dates for Ibn al-Khaww\={a}mm death in three separate places, around 934H/1528. Al-Birjandi wrote marginal commentaries on al-Faw\={a}id al-Bah\={a}iya. Tihr\={a}n\={i} notes this reference to Hajji Khalifa (\cite{Tihrani2009}, \cite{Hajji1941}).

	\item Ibn al-Khaww\={a}m himself hints, at the end of the fourth book of al-Faw\={a}id al-Bah\={a}iya, that he may have written his own commentary on this treatise elsewhere, but this commentary has not reached us \cite{Abdeljouad1986}.
	
	\item Yahya Ahmad al-K\={a}shani\cite{Kashani1320} (d.744H/1344)
	commented on Ibn al-Khaww\={a}m\={a} book in Id al-Maqid f Faw\={a}id al-Faw\={a}id
	. Abdeljouad \& Al also pointed out the frequent confusion between Ya\={h}y\={a} Im\={a}ad al-K\={a}sh\={i} and Jamsh\={i}d al-K\={a}sh\={i}, who lived at the end of the fourteenth and beginning of the fifteenth century \cite{Abdeljouad1986} and is famous in the Western world for the theorem that bears his name in trigonometry:$c^{2} = a^{2}+b^{2}-2ab \cos\theta$
	
	\item Both Yahya Ahmad al-K\={a}sh\={i} and Jamsh\={i}d al-K\={a}sh\={i} published books on arithmetic with very similar titles. The book by Yahya Ahmad al-K\={a}sh\={i} is entitled Kit\={a}b al-Hisab (Book of Arithmetic), whereas that of Jamsh\={i}d al-K\={a}sh\={i} is entitled "Talkhis Miftah Fi 'Ilam al His\={a}b" (Summary of the Key in the Science of Arithmetic). The latter seems to be the work dealing with Ibn al-Khaww\={a}m problems, and we have obtained a manuscript copy of it.
	
	\item The mathematician Turk Ram\={a}ad\={a}an - the historian $\ddot{O}$$\breve{g}r$ $\ddot{U}$yesi H$\ddot{u}$seyin ALi notes that he was still alive in 1092H/1681, after having completed his work Solution of the Sketch for the Governors in the Explanation of the Sketch of Beha-eddin Amili \cite{Hurayra1681a, Hurayra1681b}. No information is given about his date of birth or death \cite{Ogr2019}. This book has not yet been edited and remains available only in manuscript form.
	
	\item Abud El-Rahim Ibn Abu B\={a}akar Ibn Soulym\={a}n El M\={a}arachii (1068H/1658) with the title of his manuscript: Rissala al-Bah\={a}'iya Fi Al-Hissab(\cite{Maarachii}, \cite{Kahaala1957}).
	\item The Kurdish mathematician Omar ibn A\={a}mad al-M El Djali El kurdi\cite{Maayi1684}, author of the book entitled Commentaries on Problematic Places and Annotations on Symbols Resulting from Problematic Investigations.

	\item The mathematician al-Amili\={i} (953H/1547-1031H/1622) reproduced extracts from al-Faw\={a}id al-Bah\={a}iya; his seven unsolved problems appear among the thirty-three open problems of Ibn al-Khaww\={a}m (specifically problems $4, 8, 17, 18, 19, 24$, and $32$) \cite{Amili1981}.
	\item Mahdi Abdeljaouad and Hmida Hedfi commented on Ibn al-Khaww\={a}m treatise in their article Vers une \'etude des aspects historiques et math\'ematiques des probl\`emes ouverts d'Ibn al-Khaww\={a}m (XIIIe s.) \cite{Abdeljouad1986}. Despite its $23$ pages, their work remains to this day the most comprehensive attempt compared with earlier efforts by other mathematicians. They review multiple solutions, relying on the work of other scholars who had already treated some of these problems, while also identifying problems that had not been mentioned and that do not require sophisticated methods problems that today can be considered elementary. Some of these problems are even connected with questions that were later known as Fermat Last Theorem (solved by Wiles \cite{Wiles1995}), which can now be approached with elementary methods and whose demonstrations they partially reconstruct.
	
	\item The French mathematician Marre translated al-Amili book Kh\'elasat al His\={a}b  (Essence of Calculation) into French, and on this occasion reproduced problems $4, 8, 17, 18, 19, 24$, and $32$ from al-Faw\={a}id al-Bah\={a}iya as included in al-Amili  book \cite{Marre1846}. This was preceded by a german translation of Kh\'elasat al His\={a}b  by the German mathematician Nesselmann \cite{Behaeddin1843}.
	
	\item The historian of science Moustafa Mawaldi devoted the third volume of his dissertation to the mathematician ?asan al-F\={a}ris\={i} \cite{Mawaldi1989}, as well as to the thirty-three problems of Ibn al-Khaww\={a}m, treating them in great detail. Professor Mahdi Abdeljaouad kindly provided me with a copy of this doctoral dissertation dedicated to Ibn al-Khaww\={a}m\={a} list \cite{Derouiche2023}.
	
	\item The Turkish historian ihsan Fazlioglu\cite{Fazlioglu1995} has also worked on the thirty-three problems, contributing solutions inspired by the work of Mahdi Abdeljaouad and Hmida Hedfi, and adding new insights, including a discussion of Fermat conjecture.
	
	\item The mathematician and historian of mathematics Ahmed Djebbar, drawing on the article by Mahdi Abdeljaouad and Hmida Hedfi \cite{Abdeljouad1986}, reproduces in the appendix of his work Une histoire scientifique des pays d'islam : entretiens avec Jean Rosmorduc \cite{Djebbar2001} the list of the thirty-three problems of Ibn al-Khaww\={a}m, without adding commentary of his own \cite{Derouiche2022}.
\end{itemize}

\section{d'Al-Faw\={a}'id al-Bah\={a}'iya fi qaw\={a}'id al-His\={a}biya}
\begin{figure}
	\centering
	\includegraphics[width=0.6\textwidth]{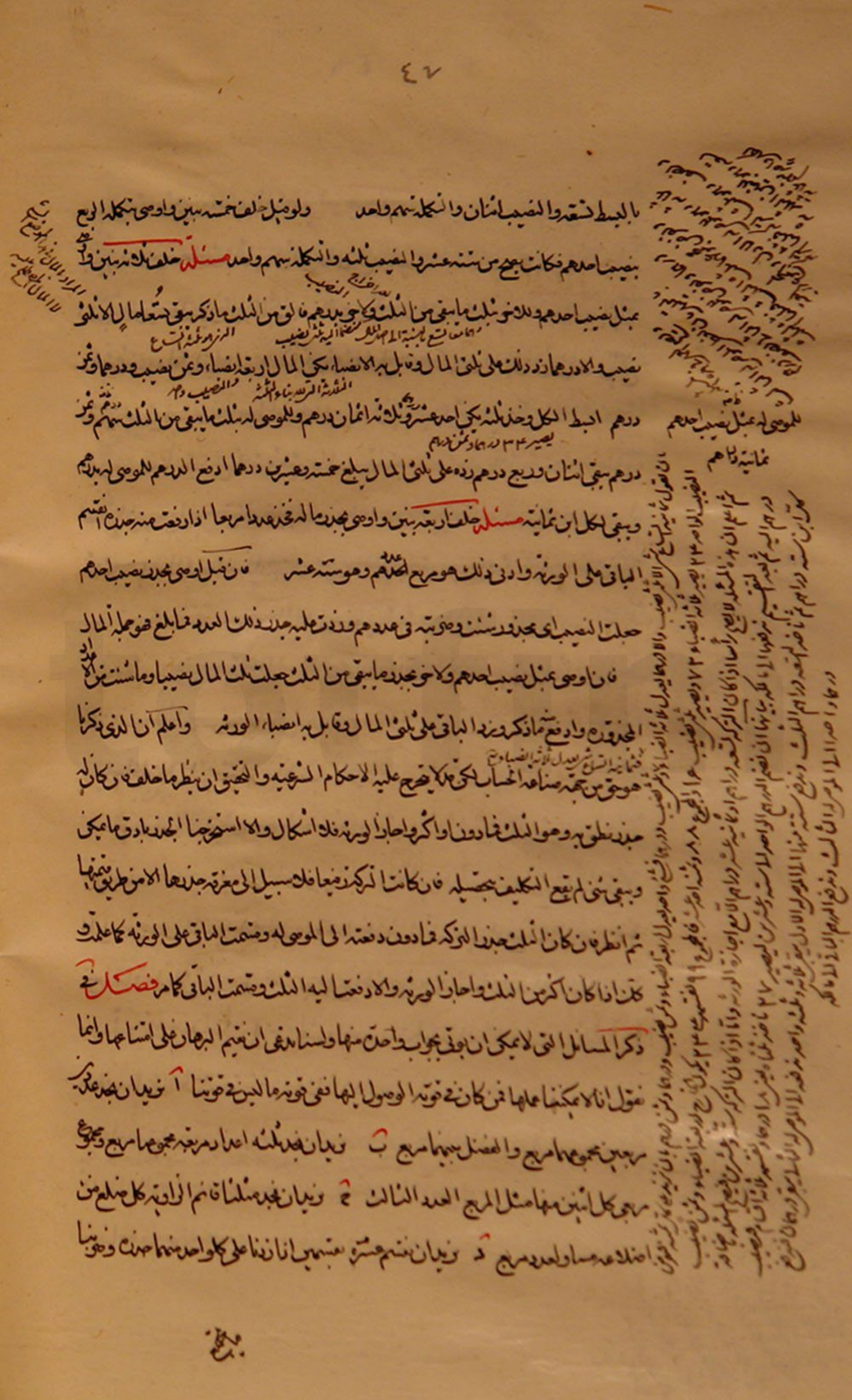}
	\caption{A snapshot of the list of the thirty-three problems from the library of Sadr Madrasa (fol. 52r).) [\cite{Khawwam}].}
	\label{figure_problems_lists}
\end{figure}

Abdeljaouad mentions five books that make up Ibn al-Khaww\={a}m treatise al-Faw\={a}'id al-Bah\={a}'iya fi qaw\={a}'id al-is\={a}'biya:

\begin{description}
	\item[Book 1 :] On multiplication, fractions, division, and root extraction.
	\item[Book 2 :] On commercial transactions and their laws.
	\item[Book 3 :] On areas and volumes.
	\item[Book 4 :] On algebra (al-jabr wal-muqbala).
	\item[Book 5 :] On the resolution of algebraic problems.
\end{description}
it is this fifth book of al-Faw\={a}'id al-Bah\={a}'iya fi qaw\={a}'id al-is\={a}biya, devoted to the resolution of algebraic problems, that will be the subject of the present section (see Figure~\ref{figure_problems_lists}).

\section{List of Ibn al-Khaww\={a}m Open Problems}

\subsection{Problems related to congruence}
\begin{remark}
	Problems $17$ and $19$ were posed by al-'\={A}mil\={i}.\cite{Abdeljouad1986}
\end{remark}
\begin{description}
	\item[Problem 1] Find two square numbers such that both their sum and their
	difference are also squares.
	\[
	\textbf{Modern formulation}: x^{2}+y^{2} =0 and x^{2}-y^{2} = 0
	\]
	
	\item[Problem 18] Find a number that admits a square root such that, if you add ten dirhams to 
	it, the result admits a square root, and if you subtract ten dirhams, the result also admits a 
	square root.
	\\
	
	\textbf{Modern formulation}: $x^{2}+10 =0$ and $x^{2}-10 =0$
	\\
	
	\item[Problem 19] A square that admits a root: if you add its root and two dirhams to it, the result is a square; if you subtract its root and two dirhams, it is again a square.
	\\
	
	\textbf{Modern formulation}: $x^{2}+\left(x+2\right) =0$ and $x^{2}-\left(x+2\right)=0$

\end{description}
\subsection{Equations of degree Four or Higher}
\begin{description}
	\item[Problem 16.] Find a square number such that, if ten times its square root and ten dirhams are subtracted from it, the result is still a perfect square. 
	\\
	\textbf{Modern formulation}: $\left(10x + 10 \right) - x^{2} =  0$.

	\item[Problem 26.] Divide a given number into two parts, one being a cube and the other a square.
	\\
	
	\textbf{Modern formulation}: $x^{3}+y^{2} = n$, given $n$
	\\
	
	\item[Problem 27.] Divide a given number into two cubic numbers.
	\\
	
	\textbf{Modern formulation}: $x^{3}+y^{3} = n$, given $n$
\end{description}
\subsection{Equations of degree Three}
Problems $15$, $21$, $22$, $25$, $28$, $29$, $32$, and $33$ fall into the category of univariate 
polynomial equations of degree four or higher. These can be transformed into Diophantine 
equations or treated via Galois theory, provided that the Galois group is solvable - which we shall 
do whenever possible.

\begin{description}
	\item[Problem 15.] Divide ten into two numbers such that, if the smaller is divided by the larger, the result is added to the latter, and the whole is multiplied by the smaller, one obtains a given number.
	\\
	
	\textbf{Modern formulation}: $x+y= 10; y\left( \frac{y}{x} + x \right) = n$, $x > y$, given $n$.
	\\
	
	\item[Problem 21.]  Find a square number such that, when its root is subtracted from it and the remainder is multiplied by its root, the result is the original square..
	\\
	
	\textbf{Modern formulation}: $\left(x^{2} - x \right)\sqrt{x^{2}-x} = x^{2}$.

	\begin{equation}
	\left( \left(x^{2}-x\right) \sqrt{x^{2}-x}\right)^{2} =  \left(x^{2}\right)^{2}
	\end{equation}
	After simplification, we obtain the equation $x^{6}-3x^{5}+2x^{4}-x^{3}$, Abdeljouad \cite{Abdeljouad1986} Proceeds by a change of variable by setting $X=x-1$ solves this equation by considering $x$ as real, and shows the impossibility of finding rational solutions.
   \\
   
	\item[Problem 22.]  Find two numbers whose difference is ten times the root of the smaller, and such that if each is divided by the other and the results are added, the sum is equal to the root of the smaller number.
	\\
	
	\textbf{Modern formulation}:  $x - y =10\sqrt{y}; \frac{x}{y}+\frac{y}{x} = \sqrt{y}$.
	\\
	
	\item[Problem 25.] Divide a given number into two parts such that, if the number is divided by each part and by their difference, the sum of all the quotients is a given number.
	\\
	
	\textbf{Modern formulation}: $x+y = a; \frac{a}{x} + \frac{a}{y} + \frac{a}{x-y} = b ; x > y $.
	\\
	
	\item [Problem 28.] A square from which its third and its root are subtracted: ten times the root of the result is equal to one third of the square plus its root.
	\\
	
	\textbf{Modern formulation}: $10\sqrt{x^{2}-(\frac{x^{2}}{3}+x) } = \frac{x^{2}}{3}+ x$.
	\\
	
	\textbf{Resolution.} This problem can be reduced by squaring both sides
%	\begin{equation}
%		\begin{split}
%			\left( 10\sqrt{x^{2}-\left(\frac{x^{2}}{3}+x  \right)} \right)^{2}= \left(\frac{x^{2}}{3}+ x\right)^{2} \\
%			\frac{200 x^{2}}{3} - 100 x = \frac{x^{4}}{9} + \frac{2 x^{3}}{3} + x^{2}
%		\end{split}
%	\end{equation}
By squaring both terms, the problem reduces to a quartic equation. This quartic can be solved by radicals using Ferrari's method, yielding the (complex) solutions. Abdeljaouad further reduces the problem to a cubic equation of the form
\[ 
X^{3}- 603X+2098 = 0
 \] 
 obtained by the change of variable $X = x+2$ \cite{Abdeljouad1986}. This cubic can be solved using Cardano's method. In both Ferrari's and Cardano's approaches, the discriminant is negative, implying that all solutions are complex. Therefore, the equation admits no rational solutions.
 \\
 
	\item [Problem 29.] Find the number of garments whose total price is a given number of dirhams, knowing that the price of a single garment plus the square root of the number of garments is equal to a given number.
	\\
	
	\textbf{Modern formulation}: $xy = a; x +\sqrt{y} = b$, given $a$ and  $b$.
	\\
	
	\item [Problem 32.]  Divide ten into two numbers such that, when each is divided by the other and the quotients are added together, the result is equal to one of the two numbers.
	\\
	
	\textbf{Modern formulation}: $x + y = 10; \frac{x}{y} + \frac{y}{x} = x$.
	\\
	
	\item [Problem 33.]  Divide ten into two numbers such that, if the larger is divided by the smaller and the quotient is multiplied by the larger, and the smaller is divided by the larger and the quotient is multiplied by the smaller, then the difference between the larger result and the smaller result is a fixed number.
	\\
	
	\textbf{Modern formulation}:   $x + y = 10;x\left( \frac{x}{y} \right)- y \left(\frac{y}{x} \right)= n$, given $n$. 
\end{description}

\subsection{Equations of degree Four or Higher }
Problems $3$, $4$, $5$, $6$, $11$, $17$, and $20$ are classified as belonging to the category of 
single-variable polynomial equations of degree $4$ or higher, which can be transformed into 
Diophantine equations or studied through Galois theory if the equation's group is 
solvable-something we will do when possible and mention whenever the occasion allows.

\begin{description}
	\item[Probleme 4.] Divide ten into two numbers such that, if you add the square root to each number and multiply the resulting sums together, you obtain a given number.
	\\
	
	\textbf{Modern formulation}: $x + y = 10; \left( x + \sqrt{x}\right)\left(y +\sqrt{y} \right)= n$, given $n$.
	\\
	
	\textbf{Resolution.} It is a matter of finding the appropriate change of variables by setting $X^{2}=x$ and $Y^{2}=y$, which transforms the problem into a system.(~\ref{eq:eq_1})
	\begin{equation}\label{eq:eq_1}
	\sigma(s,i) = \left\{
	\begin{array}{ll}
		X^{2} + Y^{2} = 10 \\
		\left(X^{2}+X\right) \left(Y^{2}+Y\right)=n
	\end{array}
	\right.
	\end{equation}
	After some algebraic manipulations between the two equations, $Y^{2} = 10-X^{2}$, as we expand in the equation: $\left(X^{2}+X\right) \left(10-X^{2}+\sqrt{10-X^{2}}\right)=n$
	\[ 
	\left(X^{2}+X\right)\sqrt{10-X^{2}} = n-\left(X^{2}+X\right)10+\left(X^{2}+X\right)X^{2}
	\]
	\item[Probleme 5.] Divide ten into two numbers such that, if each is divided by the other, the resulting quotients are added together, the sum is squared, and then multiplied by one of the two numbers, the result is a given number.
	\\
	
	\textbf{Modern formulation}: $x + y = 10; \left( \frac{x}{y} + \frac{y}{x}\right)^{2} x= n$, given $n$, we expand the expression $\left( \frac{x}{y} + \frac{y}{x}\right)^{2} x= n$, by simplifying, 
	\[ 
	\left(\frac{x^{2}+y^{2}}{yx}\right)^{2}x = n
	\]
	\[ 
	x^{4}+2x^{2}y^{2}+y^{4} = ny^{2}x 
	\]
	we substitute $y$ by $10-x$ in the equation, we obtain
	\[ 
	x^{4}+2x^{2}\left(10-x\right)^{2}+\left(10-x\right)^{4} = n\left(10-x\right)^{2}x
	\]
	\item[Problem 6.]  Divide ten into two numbers such that, if each is multiplied by its square root and the results are added together, the sum is a given number.
	\\
	
	\textbf{Modern formulation}: $x + y = 10; x\sqrt{x}+y\sqrt{y}=n$, given $n$ .
	\\
	
	\textbf{Resolution.} The problem can be written as the following system:
	\begin{equation}\label{eq:W_equivariance}
	\sigma(s,i) = \left\{
	\begin{array}{ll}
		x + y = 10 \\
		x\sqrt{x}+y\sqrt{y}=n
	\end{array}
	\right.
	\end{equation}
	 A change of variables, proposed by Abdeljaouad, consists in setting $X^{2} =x $  and $Y^{2} = y$;Under this substitution, the system (\ref{eq:W_equivariance}) becomes $X^{3}+Y^{3}=n$ and $X^{2}+Y^{2}=10$, Thus, the problem reduces to decomposing a given number $n$ 
	 as a sum of two cubes $X^{2}+Y^{2}=10$, It then suffices to perform algebraic manipulations to determine the admissible solutions
	$Y^{2} = 10-X^{2}$ for the value $Y$ of the equation $X^{3}+Y^{3}=n$,we need to expand
	\[ 
	X^{3}+\left(10-X^{2}\right)^{3} = n \\
	\]
	the equation obtained is
	\begin{equation}\label{eq:six_equ}
	X^{6} - 2 X^{3} n + X^{2} + n^{2} - 10
	\end{equation}
	Our goal is to find rational solutions of \ref{eq:six_equ}; the question is how to solve it?.
	\\
	
	\item[Problem 11.] Find two squares whose sum is ten, and such that if you add the square root to each, the results are also perfect squares.
	\\
	
	\textbf{Modern formulation}: $x^{2}+y^{2}+10$, $x+x^{2} =0; y+y^{2} = 0$.
	\\
	
	\item[Problem 17.] A person receives an inheritance of $10$ dirhams decreased by the square root of the inheritance of a second person. The inheritance of the latter is equal to $5$ dirhams decreased by the square root of the inheritance of the first person.
	\\
	
	\textbf{Modern formulation}: $y = 10 - \sqrt{y}; x = 5 - \sqrt{y}$.
	\\
	
	\textbf{Resolution.} By applying the same change of variables as in problems  $Num$, or proceeding differently if necessary, we set, $x-5 = \sqrt(5)$ and $y=10-\sqrt(10)$, By substitution, we obtain:
	\[ 
	\left(x-5\right)^{2} = \left(10-\sqrt(x)\right)
	\]
	\[
	\left(x-5\right)^{2}
	\]
 \item[Problem 20.] Find a square number which, when multiplied by itself and increased by ten times its square root and ten dirhams, becomes a number that has a square root.
 \\
 
 \textbf{Modern formulation}: $\left(x^{2} \right)^{2} + 10x + 10 = 0$.
 \\
  
\end{description}

\subsection{Diophantine equations}
A very interesting book in which recent solutions are recorded\cite{Guy2004} 
\begin{description}
	\item[Problem  2.] Find three square numbers whose sum is a square, such that the sum of the squares of any two of them is equal to the square of the third.
	\\
	
	\textbf{Modern formulation}: $x^{2} + y^{2} + z^{2} =  0; \left( x^{2}\right)^{2} + \left( y^{2}\right)^{2} = \left( z^{2}\right)^{2}; \left( x^{2}\right)^{2}+\left( 
	y^{2}\right)^{2}=\left(z^{2} \right)^{2};\left(y^{2}\right)^{2}+\left(z^{2} \right)^{2} = \left(x^{2} \right)^{2}$. 
	\\
	
	\item[Problem 3.] Find a right triangle whose sides are perfect squares.
	\\
	
	\textbf{Modern formulation}: $\left( x^{2}\right)^{2}+\left( y^{2}\right)^{2}=\left( z^{2}\right)^{2}$. 
	\\
	
	\item[Problem 8.] Find three proportional square numbers whose sum is also a square.
	\\
	
	\textbf{Modern formulation}: $\frac{x^{2}}{y^{2}} = \frac{y^{2}}{z^{2}}; x^{2}+y^{2}+z^{2}= 0$.
	\\
	
	\item[Problem 9.] Find two squares such that the product of one of them by their sum is equal to the square of the other.
	\\
	
	\textbf{Modern formulation}: $\left( x^{2}+y^{2} \right)x^{2} = \left(y^{2} \right)^{2}$.
	\\
	
	\item[Problem 12.] Find a cube such that, when one dirham is added to it, the result is another perfect cube.
	\\
	
	\textbf{Modern formulation} :  $x^{3}+1=0$ and $x^{3}+y^{3} =0$
	\\
	
	\item[Problem 13] Find two cubes such that, if each is divided by the other and the results are added together, the sum is a perfect square.
	\\
	
	\textbf{Modern formulation}: $\frac{x^{3}}{y^{3}} + \frac{y^{3}}{x^{3}} =  0 = a^{2}$.
	\\
	
	\item[Problem 14.] Find a square number such that, if each is divided by the other, the results are added together and multiplied by one of them, the outcome is a perfect square.
	\\
	
	\textbf{Modern formulation}: $x^{2}\left( \frac{x^{2}}{y^{2}} + \frac{y^{2}}{x^{2}}\right) = 0$.
	\\
	
	\item[Problem 20] Find a square number which, when multiplied by itself and increased by ten times its square root and ten dirhams, becomes a number that has a square root.
	\\
	
	\textbf{Modern formulation}: $\left(x^{2}\right)^{2}+10x+10$
	\\
	
	\item[Problem 23.] Find a cube whose difference from its square is a perfect square.
	\\
	
	\textbf{Modern formulation} $\left(x^{3}\right)-\left(x^{3}\right)^{2}=0$ and $\left(x^{3}\right)^{2}-x^{3}=0$
\end{description}
\subsection{Multiplicative problems}
\begin{description}
		\item[Problem 10] Find four square numbers such that if the first is multiplied by itself, then by the second, then by the third, and the total by the fourth, the result yields the four numbers again.
		\\
		
		\textbf{Modern formulation}: $\left(x^{2}x^{2}; x^{2}x^{2} y^{2}; x^{2}x^{2}y^{2}z^{2}; x^{2}x^{2}y^{2}z^{2}t^{2}\right) = \left(x^{2}; y^{2}; z^{2}; t^{2}\right)$. 
		\\
		
	\item [Problem  30]  Find three cubic numbers such that if the first is multiplied by the second, and the result by the third, the final product is a given number.
	\\
	
	\textbf{Modern formulation}:  $x^{3}y^{3}z^{3} = a$ ; given $a$.
	\\
	
	\item[Problem 31] An employee, whose monthly salary is a fixed but unknown number of dirhams, worked during a month for a number of days equal to the square root of the salary and received a predetermined amount of money.
	\\
	
	\textbf{Modern formulation}: $\frac{x}{30}\sqrt{x} = a$;  given $a$.

\end{description}

%%%%%%%%%%%%%%%%%%%%%%%%%%%%%%%%%%%

\section*{Acknowledgement}

The authors wish to thank Djebbar Ahmed and Mahdi Abdeljaouad for their insightful exchanges and constructive suggestions, which have led to a significant enhancement of this work.

\printbibliography

\end{document}